\documentclass[11pt]{article}

\usepackage{latexsym}       
\usepackage{amsmath} 
\usepackage{amssymb}
\usepackage{amsthm} 
\usepackage{epsfig} 

\numberwithin{equation}{section}

\newtheorem{theorem}{Theorem}[section]

\theoremstyle{definition}  

\newcommand{\p}{\partial}
\newcommand{\pd}[2]{\frac {\p #1}{\p #2}}
\newcommand{\eqnref}[1]{(\ref {#1})}

\def\nm{\noalign{\medskip}}
\newcommand{\ds}{\displaystyle}

\def\Xint#1{\mathchoice {\XXint\displaystyle\textstyle{#1}}%
{\XXint\textstyle\scriptstyle{#1}}%
{\XXint\scriptstyle\scriptscriptstyle{#1}}%
{\XXint\scriptscriptstyle\scriptscriptstyle{#1}}%
\!\int} \def\XXint#1#2#3{{\setbox0=\hbox{$#1{#2#3}{\int}$ } \vcenter{\hbox{$#2#3$ }}\kern-.6\wd0}}

\def\ddashint{\Xint\times}  

\newcommand{\Cbb}{\mathbb{C}}

\newcommand{\Rbb}{\mathbb{R}}

\newcommand{\Acal}{\mathcal{A}}
\newcommand{\Ccal}{\mathcal{C}}
\newcommand{\Ecal}{\mathcal{E}}
\newcommand{\Hcal}{\mathcal{H}}

\newcommand{\Kcal}{\mathcal{K}}
\newcommand{\Dcal}{\mathcal{D}}

\def\Be{{\bf e}}
\def\Bf{{\bf f}}
\def\Bg{{\bf g}}

\def\Bn{{\bf n}}

\def\Bu{{\bf u}}
\def\Bv{{\bf v}}
\def\Bw{{\bf w}}
\def\Bx{{\bf x}}
\def\By{{\bf y}}
\def\Bz{{\bf z}}

\def\BG{{\bf G}}

\def\BI{{\bf I}}

\def\BN{{\bf N}}

\newcommand{\Ga}{\alpha}

\newcommand{\Gd}{\delta}
\newcommand{\Ge}{\epsilon}
\newcommand{\Gve}{\varepsilon}

\newcommand{\Gvf}{\varphi}
\newcommand{\Gg}{\gamma}

\newcommand{\Gs}{\sigma}

\newcommand{\GO}{\Omega}

\newcommand{\beq}{\begin{equation}}
\newcommand{\eeq}{\end{equation}}

\begin{document}

\title{Boundary perturbations due to the presence of small linear cracks in an elastic body\thanks{\footnotesize This work was
supported by the ERC Advanced Grant Project MULTIMOD--267184 and
NRF grants No. 2009-0090250, 2010-0004091, and 2010-0017532.}}

\author{Habib Ammari\thanks{\footnotesize Department of Mathematics and Applications, Ecole Normale Sup\'erieure,
45 Rue d'Ulm, 75005 Paris, France (habib.ammari@ens.fr).} \and
Hyeonbae Kang\thanks{Department of Mathematics, Inha University, Incheon
402-751, Korea (hbkang@inha.ac.kr, hdlee@inha.ac.kr, jisunlim@inha.ac.kr).}  \and
Hyundae Lee\footnotemark[3]  \and Jisun Lim\footnotemark[3]}

\maketitle

\begin{abstract}
In this paper, Neumann cracks in elastic bodies are
considered. We establish a rigorous asymptotic expansion for the boundary
perturbations of the displacement (and traction) vectors that are due to the presence of a small elastic linear
crack. The formula reveals that the leading order term is $\Gve^2$ where $\Gve$ is the length of the crack, and the $\Gve^3$-term vanishes. We obtain an asymptotic expansion of the elastic potential energy as an immediate consequence of the boundary perturbation formula. The derivation is based on layer potential
techniques. It is expected that the formula would lead to very
effective direct approaches for locating a collection of small
elastic cracks and estimating their sizes and orientations.
\end{abstract}


\noindent {\footnotesize {\bf Mathematics subject classification
(MSC2000)}: 35B30, 74B05}


\section{Introduction}

The displacement (or traction) vector can be perturbed due to the presence of a small crack in an elastic
medium. The aim of this paper is to derive an asymptotic formula for
the boundary perturbations of the displacement as the length of the crack tends to zero. The focus is on cracks with homogeneous Neumann boundary
conditions, {\it i.e.},  perfectly insulating cracks. We
consider the linear isotropic elasticity system in two dimensions and assume that the crack is a
line segment of small size. The derivation of the asymptotic formula is based on layer potential techniques.

The paper extends recent asymptotic results that have been used
for an efficient imaging of small defects. In \cite{AKLP} an
electrostatic model, where the crack is perfectly conducting, was considered and an asymptotic expansion of
the boundary perturbations that are due to the presence of a small
linear crack was derived. The asymptotic formula leads us to efficient algorithms to detect cracks using
 boundary measurements \cite{AGKPS, AKLP}. Their resolution and stability of the algorithms with respect to medium
and measurement noises were investigated in \cite{AGJK, AGS}.
There were some work on boundary perturbation due to the presence of
small inclusions in linear elasticity; the effect of small
inclusions on boundary measurements has been studied in
\cite{AK2004, AKNT}. The effect of thin elastic inclusions on
boundary measurements was quantified in \cite{BF,BFV}. Direct
reconstruction algorithms for locating small or thin elastic
defects were developed in \cite{ACI, BFKL, KKL1, KKL2}. We
emphasize that the results of this paper (on cracks) can {\it not}
be obtained as a limiting case of thin inclusions.

The results of this paper reveals that the leading order term of the boundary perturbation
is $\Gve^2$ where $\Gve$ is the length of the crack, and its intensity is given by the traction
force of the background solution on the crack (see Theorem \ref{thm-neuman}.) We also prove that
the $\Gve^3$-order term vanishes. By integrating the boundary perturbation
formula against the given traction, we are able to derive an asymptotic expansion for
the perturbation of the elastic potential energy, which is an improvement over already existing
results \cite{GN, NDTP} (see the discussion at the end of Section 5).

The boundary perturbation formula derived in this paper carries information about the location,
size, and orientation of the crack, and we expect,  as in the electrostatic case, that the formula will provide a
powerful tool to solve the inverse problem of identifying the
cracks in terms of boundary measurements. The implementation of imaging
algorithms based on the present expansion and the analysis of
their resolution and stability will be the subject of a
forthcoming paper.

The paper is organized as follows. In section \ref{sect2}, a
representation formula for the solution of the problem in the
presence of a Neumann crack is derived. Section \ref{sect3} is
devoted to making explicit the hyper-singular character involved in the
representation formula. Using analytical results for the finite
Hilbert transform, we derive in section \ref{sect4} an asymptotic
expansion of the effect of a small Neumann crack on the boundary
values of the solution. Section \ref{sect5} aims to derive the
topological derivative of the elastic potential energy functional.
Appendix contains technical calculation of the double layer potential.







\section{A representation formula} \label{sect2}

Let $\Omega \subset \Rbb^2$ be an open bounded domain, whose
boundary $\p \Omega$ is of class $\Ccal^{1, \Ga}$ for some $\Ga>0$. We assume that $\GO$
is a homogeneous isotropic elastic body so that its elasticity
tensor $\Cbb=(C_{ijkl})$ is given by
\begin{equation}
C_{ijkl}= \lambda \delta_{ij}\delta_{kl} + \mu (\delta_{ik}\delta_{jl}+ \delta_{il}\delta_{jk})
\end{equation}
with the Lam\'{e} coefficients $\lambda$ and $\mu$  satisfying
$\mu >0$ and $\lambda + \mu >0$. Let ${\gamma}_{\Gve} \subset
\Omega$ be a small straight crack with  size $\Gve$, located at
some fixed distance $d_0$ from $\partial \Omega$,  {\it i.e.},
$$\mbox{dist}({\gamma}_{\Gve},
\partial \Omega) \ge d_0.$$
We denote by $\Be^\perp$ a unit normal to $\Gg_\Gve$.

Let
$$\Psi := \bigg\{ \psi: \partial_i \psi_j + \partial_j \psi_i = 0, \quad
1\leq i,j\leq 2 \bigg\},$$ or equivalently,
$$\Psi = \mbox{span}\bigg\{\begin{bmatrix} 1 \\
0
\end{bmatrix}, \begin{bmatrix} 0
\\ 1 \end{bmatrix}, \begin{bmatrix} y \\ -x \end{bmatrix}
\bigg\}.$$ Introduce the space
$$
L^2_\Psi(\partial \Omega) := \bigg\{ \mathbf{f} \in L^2(\partial
\Omega):  \int_{\partial \Omega} \mathbf{f} \cdot {\psi} \,
d\sigma = 0 \quad \mbox{for all } \psi \in \Psi \bigg\}.
$$
Let $\Bu_\Gve$ be the displacement vector caused by the traction
$\Bg \in L^2_\Psi(\partial \Omega)$ applied on the boundary $\p\GO$ in the presence of $\Gg_\Gve$.
Then $\Bu_\Gve$ is the solution to
\begin{equation}\label{prob1}
\begin{cases}
\nabla \cdot \sigma(\Bu_\Gve) = 0 \quad  & \mbox{in}\; \Omega \setminus \bar{\gamma}_{\Gve}, \\
\sigma(\Bu_\Gve)\,\Bn = \Bg   & \mbox{on} \; \partial \Omega,\\
\sigma(\Bu_\Gve)\,\Be^\perp = 0   & \mbox{on} \; \gamma_{\Gve},
\end{cases}
\end{equation}
where $\Bn$ is the outward unit  normal to $\p\GO$ and
$\Gs(\Bu_\Gve)$ is the stress defined by
\begin{equation}
\sigma(\Bu_\Gve) = \Cbb \nabla^s \Bu_\Gve := \frac{1}{2} \Cbb (\nabla \Bu_\Gve +\nabla \Bu_\Gve^T ) .
\end{equation}
Here, $\nabla^s \Bu_\Gve=\frac{1}{2}  (\nabla \Bu_\Gve +\nabla \Bu_\Gve^T )$ is the
strain tensor and the superscript $T$ denotes the transpose of a
matrix. Note that the functions in $\Psi$ are solutions to the
homogeneous problem (\ref{prob1}) with $\Bg=0$. So we impose
orthogonality condition on
 $\Bu_\Gve$ to guarantee the uniqueness of a solution to \eqnref{prob1}:
\begin{equation}\label{nor_cond}
\int_{\p\GO} \Bu_\Gve \cdot \psi \, d\Gs=0 \quad
\mbox{for all } \psi \in \Psi .
\end{equation}
Let  $\Bu_0$ be the solution in the absence of the crack, {\it
i.e.}, the solution to
\begin{equation}\label{modelu0}
\begin{cases}
\nabla \cdot \sigma(\Bu_0) = 0 \quad  & \mbox{in}\; \Omega , \\
\sigma(\Bu_0)\,\Bn = \Bg   & \mbox{on} \; \partial \Omega,
\end{cases}
\end{equation}
with the orthogonality condition: $\Bu_0 |_{\partial \Omega} \in
L^2_\Psi(\partial \Omega)$ (or equivalently, \eqnref{nor_cond}
with $\Bu_\Gve$ replaced with $\Bu_0$).

It is well-known that the solution $\Bu_\Gve$ to \eqnref{prob1}
belongs to $H^1(\GO \setminus \Gg_\Gve)$. In fact, we have
 \beq\label{bound}
 \| \Bu_\Gve \|_{H^1(\GO \setminus \Gg_\Gve)} \le C
 \eeq
for some $C$ independent of $\Gve$. To see this, we introduce the potential energy functional
 \beq\label{energy}
 J_\Gve[\Bu] := -\frac{1}{2} \int_{\GO \setminus \Gg_\Gve} \Gs(\Bu) : \nabla^s
 \Bu.
 \eeq
The solution $\Bu_\Gve$ of (\ref{prob1}) is the maximizer
of $J_\Gve$, {\it i.e.},
 \beq
 J_\Gve[\Bu_\Gve] = \max J_\Gve[\Bu],
 \eeq
where the maximum is taken over all $\Bu \in H^1(\GO \setminus
\Gg_\Gve)$ satisfying $\Gs(\Bu) \,\Bn = \Bg$ on $\p\GO$ and
$\Gs(\Bu)\,\Be^\perp = 0$ on $\Gg_\Gve$.  Let $\Bv$ be a smooth
function with a compact support in $\GO$ such that
$\Gs(\Bv)\,\Be^\perp = - \Gs(\Bu_0)\,\Be^\perp$ on $\Gg_\Gve$. We may
choose  $\Bv$ so that $J_\Gve[\Bv]$ is independent of $\Gve$.
Since $0 \ge J_\Gve[\Bu_\Gve] \ge J_\Gve[\Bu_0 + \Bv]$, we have
 $$
 \| \nabla^s \Bu_\Gve \|_{L^2(\GO \setminus \Gg_\Gve)} \le C.
 $$
We then have from the Korn's inequality that there is a constant $C$ independent of $\Gve$ such that
 \beq\label{korn}
 \| \Bu_\Gve - \Bu_0 \|_{H^1(\GO \setminus \Gg_\Gve)} \le C (\| \nabla^s (\Bu_\Gve - \Bu_0) \|_{L^2(\GO\setminus \Gg_\Gve)} + \| \Bu_\Gve - \Bu_0 \|_{H^{1/2}(\p\GO)}).
 \eeq
Since $\| \Bu_\Gve - \Bu_0 \|_{H^{1/2}(\p\GO)}$ is bounded regardless of $\Gve$ as we shall show later (Theorem \ref{thm-neuman}), we obtain \eqnref{bound}.

Let
\begin{equation}\label{phi}
\varphi_\Gve(\Bx):=\Bu_\Gve|_{+}(\Bx) - \Bu_\Gve|_{-}(\Bx), \qquad \Bx \in \gamma_{\Gve},
\end{equation}
where $+\,(\mbox{resp. } -)$ indicates the limit on the crack
$\gamma_{\Gve}$ from the given normal direction $\Be^\perp$ (resp.
opposite direction), {\it i.e.},
$$ \Bu_{\pm}(x):= \lim_{t\to 0} \Bu(\Bx\pm t \Be^\perp).$$
We sometimes denote $\sigma(\Bu)\,\Bn$, the traction on $\p\Omega$ (or on $\Gg_\Gve$),  by $\p \Bu /\p \nu$, {\it i.e.},
\begin{equation}\label{conormal1}
\frac{\p \Bu}{\p \nu} := \lambda (\nabla \cdot \Bu) \Bn + \mu( \nabla \Bu + \nabla \Bu^T ) \Bn \quad \mbox{on} \; \p \Omega.
\end{equation}
If ${\Phi}=(\Phi_{ij})_{2\times 2}$ is the Kelvin matrix of the fundamental solutions of Lam\'{e} system, {\it i.e.},
\begin{equation}\label{lamefund}
   \Phi_{ij}(\Bx):=\frac{A}{2\pi}\delta_{ij}\log|\Bx|- \frac{B}{2\pi}\frac{x_ix_j}{|\Bx|^2},\qquad \Bx \neq 0 \in \mathbb{R}^2,
\end{equation}
where
\beq
A= \frac{\lambda+3\mu}{2\mu(\lambda+2\mu)} \quad\mbox{and}\quad
B= \frac{\lambda+\mu}{2\mu(\lambda+2\mu)},
\eeq
then the solution $\Bu_\Gve$ to \eqref{prob1} is represented as
\begin{align}\label{uep_integral}
\Bu_\Gve(\Bx)&=
 \int_{\p \Omega} \frac{\p  \Phi}  {\p \nu_\By} (\Bx-\By) \Bu_\Gve(\By)\,d\sigma(\By) -\int_{\p \Omega}  \Phi(\Bx-\By) \frac{\p \Bu_\Gve} {\p \nu}(\By) \,d\sigma(\By)\nonumber\\
&\quad  -   \int_{\gamma_{\Gve}} \frac{\p \Phi } {\p \nu_\By}(\Bx-\By) \varphi_\Gve(\By) \,d\sigma(\By), \qquad \Bx \in \GO\setminus \Gg_\Gve.
\end{align}
The solution $\Bu_0$ to \eqref{modelu0} is represented as
 $$
 \Bu_0(\Bx)=
 \int_{\p \Omega} \frac{\p  \Phi}  {\p \nu_\By} (\Bx-\By) \Bu_0(\By)\,d\sigma(\By) -\int_{\p \Omega}  \Phi(\Bx-\By) \frac{\p \Bu_0} {\p \nu}(\By) \,d\sigma(\By).
 $$

Let
 \beq
 \Bw_\Gve:= \Bu_\Gve - \Bu_0.
 \eeq
Since $\frac{\p \Bu_\Gve} {\p \nu}= \frac{\p \Bu_0} {\p \nu}$ on $\p\GO$, by subtracting above two identities, we have
 \beq\label{intgralform}
 \Bw_\Gve (\Bx) - \Dcal_\GO[\Bw_\Gve](\Bx) = - \Dcal_\Gve [\Gvf_\Gve](\Bx), \qquad \Bx \in \GO,
 \eeq
where the double layer potentials $\Dcal_\GO$ and $\Dcal_\Gve$ are
defined by
 \beq
 \Dcal_\GO[\Bw_\Gve](\Bx):=  \int_{\p \Omega} \frac{\p  \Phi}
  {\p \nu_\By} (\Bx-\By) \Bw_\Gve(\By)\,d\sigma(\By), \quad \Bx \in \GO
 \eeq
and
 \beq
 \Dcal_\Gve [\Gvf_\Gve](\Bx):= \int_{\gamma_{\Gve}} \frac{\p \Phi } {\p \nu_\By}(\Bx-\By) \varphi_\Gve(\By) \,d\sigma(\By), \quad \Bx \in \GO.
 \eeq

Let $\Bn_\Bx$ denote the outward normal to $\p\GO$ at $\Bx \in \p\GO$ and let
 $$
 \Dcal_\GO[\Gvf] \big|_{-}(\Bx) = \lim_{t \to 0+} \Dcal_\GO[\Gvf](\Bx+t \Bn_\Bx).
 $$
Then, it is well known (see, for example, \cite{book2}) that
 $$
 \Dcal_\GO[\Gvf] \big|_{-}(\Bx) = \left(\frac{1}{2} I + \Kcal_\GO \right) [\Gvf](\Bx), \quad\mbox{a.e. } \Bx \in \p\GO,
 $$
where $\Kcal_\GO$ is the boundary integral operator defined by
 $$
 \Kcal_\GO[\Gvf](\Bx) := \mbox{p.v.} \int_{\p \Omega} \frac{\p  \Phi}  {\p \nu_\By} (\Bx-\By) \Bw_\Gve(\By)\,d\sigma(\By), \quad \Bx \in \p\GO,
 $$
and $I$ is the identity operator. Here, p.v. stands for the Cauchy principal value. It then follows from \eqnref{intgralform} that
 \beq
 \left(-\frac{1}{2} I + \Kcal_\GO \right) [\Bw_\Gve] (\Bx) = \Dcal_\Gve [\Gvf_\Gve](\Bx), \quad \Bx \in \p\GO.
 \eeq
Since $-\frac{1}{2} I + \Kcal_\GO$ is invertible on
$L^2_\Psi(\p\GO)$ (see, for instance, \cite{DKV88}), we have
 \beq\label{14}
 \Bw_\Gve(\Bx)= \int_{\gamma_{\Gve}} \frac{\p} {\p \nu_\By} (-\frac{1}{2} I + \Kcal_\GO)^{-1}
 \left[\Phi(\cdot-\By)\right](\Bx) \varphi_\Gve(\By) \,d\sigma(\By), \quad \Bx \in \p \GO.
 \eeq
Note that
 \beq\label{fundNeumann}
 (-\frac{1}{2} I + \Kcal_\GO)^{-1} \left[\Phi(\cdot-\By)\right](\Bx) = \BN(\Bx, \By), \quad \Bx \in \p\GO, \ \ \By \in \GO,
 \eeq
modulo a function in $\Psi$, where
$\BN(\Bx, \By)$ is the Neumann function for the Lam\'{e} system on
$\GO$, namely, for $\By \in \GO$, $\BN(\Bx, \By)$ is the solution
to
\begin{equation}
\begin{cases}
 \nabla \cdot \sigma(\BN(\cdot,\By)) = - \delta_{\By} \BI \quad & \mbox{in}\; \Omega, \\
 \nm
 \ds \sigma(\BN(\cdot, \By))\,\Bn = - \frac{1}{|\p \Omega|} \BI    \quad & \mbox{on} \; \p \Omega,
\end{cases}  
\label{neumann_def}
\end{equation}
subject to the orthogonality condition:
$$\int_{\p \Omega} \BN(\Bx, \By) \cdot \psi(\Bx) \, d\Gs(\Bx) = 0 \quad \mbox{for all }
\psi \in \Psi. $$
Here, $\BI$ is the $2 \time 2$ identity matrix. See \cite{book2, AKL2009,Kup} for properties of the Neumann
function and a proof of (\ref{fundNeumann}). Thus we obtain from
\eqnref{14} that
 \beq\label{Neumann_repr}
 \Bu_\Gve (\Bx) = \Bu_0(\Bx)+ \int_{\gamma_{\Gve}} \frac{\p} {\p \nu_\By} \BN(\Bx, \By) \varphi_\Gve(\By) \,d\sigma(\By), \quad \Bx \in \p \GO.
 \eeq

We now describe the scheme to derive an asymptotic expansion of $\Bu_\Ge-\Bu_0$ on $\p\GO$. Since $\pd{\Bu_\Ge}{\nu} =\Gs(\Bu_\Ge)\Be^\perp=0$ on $\Gg_\Gve$, we use \eqnref{intgralform} to obtain
 \beq\label{inteqn1}
 \pd{\Bu_0}{\nu} + \pd{}{\nu}\Dcal_\GO[\Bw_\Gve] = \pd{}{\nu} \Dcal_\Gve [\Gvf_\Gve] \quad\mbox{on } \Gg_\Ge.
 \eeq
We solve this integral equation for $\Gvf_\Gve$ and then
substitute it into \eqnref{Neumann_repr} to derive an asymptotic
expansion of $\Bu_\Gve$ as $\Gve \to 0$.

\section{Derivation of an explicit integral equation} \label{sect3}

In view of \eqnref{inteqn1}, we need to compute $\frac{\p  }{\p \nu_\Bx}\left(\frac{\p \Phi} {\p \nu_\By}(\Bx-\By) \right)$ on $\Gg_\Gve$. As before $\Be^\perp$ is the unit normal to $\Gg_\Gve$, and denoted by $\Be^\perp=(n_1, n_2)$. It is worth mentioning that $n_1$ and $n_2$ are constant since $\Gg_\Gve$ is a line segment. It is convenient to use the following expression of the conormal derivative:
 \beq
 \frac{\p \Bu}{\p \nu} =  T(\p)\Bu ,
 \eeq
where the operator $T(\p) =T(\p_1,\p_2)$, where $\p_j=\pd{}{x_j}$, is defined by
\begin{equation}\label{op-T}
T(\xi_1,\xi_2):=\begin{bmatrix}
 (\lambda+2\mu) n_1 \xi_1 + \mu n_2 \xi_2  &
 \mu n_2 \xi_1 + \lambda n_1 \xi_2   \smallskip \\
   \lambda n_2 \xi_1 + \mu n_1 \xi_2 &
    \mu n_1 \xi_1 + (\lambda+2\mu) n_2 \xi_2
\end{bmatrix}.
\end{equation}

We first obtain the following formula whose derivation will be
given in Appendix \ref{appendix}. For $\Bx \neq \By$, we have
\begin{align}\label{repDij}
\ds \Big(\frac{\p \Phi}{\p \nu_\By}(\Bx-\By)\Big)_{ij}=&
 \left[ a\,\delta_{ij} +b\frac{(x_i-y_i)(x_j-y_j)}{|\Bx-\By|^2}\right]
 \sum_{l=1}^2 \frac{n_l(x_l-y_l)}{|\Bx-\By|^2} \nonumber \\
&\qquad   -\,a \frac{n_j(x_i-y_i)-n_i(x_j-y_j)}{|\Bx-\By|^2},
  \end{align}
where
 \beq
 a= - \frac{\mu}{2\pi( \lambda+2\mu)}, \quad
 b =  - \frac{( \lambda+\mu)}{\pi( \lambda+2\mu)} .
 \eeq
Let $v_{ij}=(\frac{\p \Phi}{\p\nu_\By}(\Bx-\By))_{ij}$ for convenience and let
\begin{equation}\label{W-fn}
W(\Bx-\By):=\frac{\p  }{\p \nu_\Bx}\left(\frac{\p \Phi} {\p \nu_\By}(\Bx-\By)\right).
\end{equation}
Then one can use \eqnref{op-T} to derive
\begin{eqnarray*}
W(\Bx-\By)_{11}
&=& n_1[ (\lambda +2\mu)  \p_1v_{11}
+\lambda {\p_2}v_{21}] + n_2 ( \mu{\p_2}v_{11}+  \mu{\p_1}v_{21}),\\
W(\Bx-\By)_{12}
&=&n_1[(\lambda+2\mu) \p_1v_{12} + \lambda{\p_ {x_2}}v_{22} ]+n_2 ( \mu\p_2v_{12} +  \mu{\p_1}v_{22}),\\
W(\Bx-\By)_{21}
&=& n_1 (\mu\p_2v_{11}+\mu\p_1v_{21})
+n_2[ \lambda \p_1v_{11}+(\lambda +2\mu)\p_2v_{21}],\\
W(\Bx-\By)_{22}
&=&n_1(\mu {\p_2}v_{12}+
\mu {\p_1}v_{22}) +n_2[\lambda {\p_1}v_{12}
+(\lambda+2\mu) {\p_2}v_{22}].
\end{eqnarray*}

Since  the crack which we consider is a line segment with length $\Gve$ in the domain $\Omega \subset \mathbb{R}^2$,
we may assume, after rotation and translation if necessary, that it is given by
\begin{equation}\label{linecrack}
\gamma_\Gve = \{ (x_1, 0) :  -\Gve/2 \le x_1 \le \Gve/2 \}.
\end{equation}
In this case, one can check that
\begin{align*}
\p_2v_{11}&= \left(a+b\frac{(x_1-y_1)^2}{|\Bx-\By|^2}\right)
\frac{1}{|\Bx-\By|^2}=\frac{a+b}{(x_1-y_1)^2},\\
\p_1v_{21}&= a\left( \frac{1}{|\Bx-\By|^2 } -2\frac{x_1-y_1}{|\Bx-\By|^2}\frac{x_1-y_1}{|x_1-y_1|^2}\right) = -\frac{a}{(x_1-y_1)^2},\\
\p_1v_{12}&= -a\left( \frac{1}{|\Bx-\By|^2 } -2\frac{x_1-y_1}{|\Bx-\By|^2}\frac{x_1-y_1}{|x_1-y_1|^2}\right) = \frac{a}{(x_1-y_1)^2},\\
\p_2v_{22}&= \frac{a}{(x_1-y_1)^2} \\
\p_1v_{11}&= \p_2v_{21} = \p_2v_{12} = \p_1v_{22} = 0 .
\end{align*}
Since $\Be^\perp=(n_1,n_2)=(0,1)$, we have
\begin{eqnarray*}
W(\Bx-\By)_{11} &=& \mu \p_2v_{11}+\mu \p_1v_{21}\\
&=& \mu(a+b) \frac{1}{(x_1-y_1)^2} - \mu a\frac{1}{(x_1-y_1)^2}
=  \frac{\mu b}{(x_1-y_1)^2},\\
W(\Bx-\By)_{12} &=& \mu \p_2v_{12}+\mu \p_1v_{22}
= 0,\\
W(\Bx-\By)_{21} &=& \lambda  \p_1v_{11}+ (\lambda+2  \mu) \p_2v_{21}= 0,\\
W(\Bx-\By)_{22} &=& \lambda  \p_1v_{12}+ (\lambda+2\mu) \p_2v_{22}\\
 &=& \lambda  a \frac{1}{(x_1-y_1)^2}+(\lambda+2\mu)a \frac{1}{(x_1-y_1)^2}
 =  \frac{2(\lambda+\mu) a}{(x_1-y_1)^2},
\end{eqnarray*}
that is,
\begin{equation}
W(\Bx-\By)=  \frac{1}{(x_1-y_1)^2}
 \begin{bmatrix}
 -\frac{\mu (\lambda+\mu)}{\pi(\lambda+2\mu)}& 0 \\  0 &  -\frac{\mu(\lambda+\mu)}{\pi(\lambda+2\mu)}
 \end{bmatrix} .
\end{equation}
Note that
\begin{equation}\label{c-lam}
\frac{\mu(\lambda+\mu)}{\lambda+2\mu}=\frac{E}{4}
\end{equation}
where $E$ is the Young's modulus in two dimensions. So, we have
\begin{equation}
W(\Bx-\By)=  -\frac{E}{4\pi}\frac{1}{(x_1-y_1)^2} \BI.
\end{equation}

So far we have shown that if $\Gg_\Gve$ is given by \eqnref{linecrack}, then
\beq\label{ddb2}
 \frac{\p } {\p \nu_\Bx}  \mathcal{D}_{\Gve}[\varphi_\Gve](\Bx)= - \frac{E}{4\pi}
 \int_{-\Gve/2}^{\Gve/2} \frac{\varphi_{\Gve}(y)}{(x-y)^2}
   \,   dy,\quad \Bx=(x,0), \ -\Gve/2<x<\Gve/2.
  \eeq
Here the integral is hyper-singular and should be understood as a finite part in the sense of Hadamard, which will be defined in the next section. So the integral equation \eqnref{inteqn1} becomes
 \beq
 \frac{1}{\pi} \int_{-\Gve/2}^{\Gve/2} \frac{\varphi_{\Gve}(y)}{(x-y)^2} \,dy= - \frac{4}{E} \Bf(x), \quad -\Gve/2<x<\Gve/2,
 \eeq
where
 \beq
 \Bf(x)=  \pd{\Bu_0}{\nu}(x,0) + \pd{}{\nu}\Dcal_\GO[\Bw_\Gve](x,0) .
 \eeq
Define
\beq
\Bf_\Gve(x) := \Bf (\frac{\Gve}{2} x),
\eeq
and
\beq\label{psiGve}
\psi_\Gve(x):= \frac{2}{\Gve} \Gvf_\Gve (\frac{\Gve}{2} x), \quad -1< x <1.
\eeq
Then the scaled integral equation is
\beq\label{scaleeqn}
 \frac{1}{\pi} \int_{-1}^{1} \frac{\psi_{\Gve}(y)}{(x-y)^2}\,dy = - \frac{4}{E} \Bf_\Gve (x), \quad -1<x<1,
 \eeq
which we solve in the next section.

\section{Asymptotic expansion} \label{sect4}
The integral in \eqnref{scaleeqn} is understood as a finite-part in the
 sense of Hadamard \cite{Mart92, Mart96}: for $\psi \in \mathcal{C}^{1,\alpha}(-1,1)  ~(0<\alpha\le 1)$
\begin{equation}\label{HM}
\ddashint_{-1}^{1} \frac{\psi(y)}{(x-y)^2}\,dy = \lim_{\Gd
\rightarrow 0} \Big[\int_{-1}^{x-\Gd} \frac{\psi(y)}{(x-y)^2}dy  +
\int_{x+\Gd}^1 \frac{\psi(y)}{(x-y)^2}dy -
\frac{2\psi(x)}{\Gd}\Big].
\end{equation}
Define
\begin{equation}\label{eq:pv}
\Acal[\psi](x): = \frac{1}{\pi} \ddashint_{-1}^{1} \frac{\psi(y) }{(x-y)^2}\,dy, \quad |x|<1.
\end{equation}
It is known (\cite{Ioak, Mart89}) that \beq\label{aandh}
\Acal[\psi](x) =- \frac{d}{dx} \Hcal[\psi](x), \eeq where $\Hcal$
is the (finite) Hilbert transform, {\it i.e.}, \beq \Hcal[\psi](x)
= \mbox{p.v.} \frac{1}{\pi}\int_{-1}^{1}  \frac{\psi(y)}{x-y}\,dy.
\eeq
More properties of finite-part integrals and principal-value integrals
can be found in \cite{Kup, Mart92, Mart96, Musk,Tric, Zabr}.

If $\psi(-1)=\psi(1)=0$, we have from \eqnref{aandh} that
\begin{equation}\label{hilbt}
\Acal[\psi](x)= -\Hcal[\psi'](x).
\end{equation}
Thus we can invert the operator $\Acal$ using the properties of $\Hcal$.
The set $\mathcal{Y}$, given by
\begin{equation} \label{xeps}
  \mathcal{Y} \,=\,
\ds \bigg\{ \varphi: \int_{-1}^1 \sqrt{1 - x^2}\,
|\varphi(x)|^2 \, dx < +\infty
  \bigg\},
  \end{equation}
 is a Hilbert space with the norm
$$ || \varphi ||_{\mathcal{Y}} = \bigg( \int_{-1}^1 \sqrt{1 - x^2}\,
|\varphi(x)|^2 \, dx  \bigg)^{1/2}.$$
It is well known (see, for example, \cite[section 5.2]{AKL2009}) that $\Hcal$ maps $\mathcal{Y}$ onto itself and its null space is the one dimensional space generated by $1/\sqrt{1-x^2}$. Therefore, if we define
\beq\label{yeps}
\begin{array}{lll}
 \mathcal{X} &=&\ds
\bigg\{  \psi \in \mathcal{C}^0 \left(\,\left[\, -1,
1 \,\right]\,\right): \psi^{\prime} \in \mathcal{
Y},~ \psi(-1)=\psi(1)=0 \bigg\},
\end{array}
\eeq
where  $\psi^{\prime}$ is the distributional derivative of $\psi$, then $\mathcal{A}: \mathcal{X}\rightarrow \mathcal{Y}$ is invertible. We note that $\mathcal{X}$ is a Banach space with the norm
\[\ds
|| \psi ||_{\mathcal{X}} \,=\,
||\psi||_{L^\infty} + ||
\psi^\prime||_{\mathcal{Y}} .
\]
Using the Hilbert inversion formula (see, for example,
\cite[section 5.2]{AKL2009}), we can check
\begin{align}
\mathcal{A}^{-1}[ 1 ](x)& =  -\sqrt{1-x^2},\label{invA1} \\
\mathcal{A}^{-1}[y ](x)& = -\frac{x}{2}\sqrt{1-x^2}.\label{invAy}
\end{align}

The equation \eqref{scaleeqn} can be written as
\beq\label{a_equation}
 \mathcal{A}[ \psi_{\Gve}](x) = -\frac{4}{E} \Bf_\Gve (x), \quad -1<x<1.
\eeq
The Taylor expansion yields
\beq
\pd{\Bu_0}{\nu}(\frac{\Gve}{2} x,0) =
\frac{\p \Bu_{0}} {\p \nu}(0) +
\frac{\Gve x}{2} \frac{\p^2  \Bu_{0}}  {\p t \p \nu}(0)
+ e_1(x),
\eeq
where ${\p}/{\p t}$ denotes the tangential derivative on $\Gg_{\Gve}$. The remainder term $e_1$ satisfies
 $$
 |e_1(x)| \le C \Gve^2 |x|^2,
 $$
and in particular,
 \beq
 \| e_1 \|_{\mathcal{Y}} \le C \Gve^2.
 \eeq
On the other hand, since
$$
\pd{}{\nu}\Dcal_\GO[\Bw_\Gve](\frac{\Gve}{2} x,0) = \int_{\p\GO}
\frac{\p^2 \BN}{\p \nu_\Bx\p \nu_\By}((\frac{\Gve}{2} x,0),\By)
\Bw_\Gve(\By) \,d\sigma(\By),
$$
and $\Gg_\Gve$ is away from $\p\GO$, one can see that \beq \left\|
\pd{}{\nu}\Dcal_\GO[\Bw_\Gve](\frac{\Gve}{2} \cdot,0)
\right\|_{\mathcal{Y}} \le C \| \Bw_\Gve \|_{L^\infty(\p\GO)}.
\eeq Therefore, we have \beq \Bf_\Gve (x)= \frac{\p \Bu_{0}} {\p
\nu}(0) + \frac{\Gve x}{2} \frac{\p^2  \Bu_{0}}  {\p t \p \nu}(0)
+ e(x), \eeq where $e$ satisfies \beq \| e \|_{\mathcal{Y}} \le C
\bigg(\Gve^2 + \| \Bw_\Gve \|_{L^\infty(\p\GO)} \bigg). \eeq

We now obtain from \eqnref{a_equation} that
\beq
\psi_\Gve (x)= -\frac{4}{E} \left[ \frac{\p \Bu_{0}} {\p \nu}(0) \Acal^{-1}[1](x) + \frac{\Gve}{2} \frac{\p^2  \Bu_{0}}  {\p t \p \nu}(0) \Acal^{-1}[y](x) + \Acal^{-1}[e](x) \right].
\eeq
Note that $\Ecal_1(x) = \Acal^{-1}[e](x)$ satisfies
$$
\| \Ecal_1 \|_{\mathcal{X}} \le C \| e \|_{\mathcal{Y}} \le C
\bigg(\Gve^2 + \| \Bw_\Gve \|_{L^\infty(\p\GO)} \bigg),
$$
and in particular,
\beq
\| \Ecal_1 \|_{L^\infty(-1,1)} \le C
\bigg(\Gve^2 + \| \Bw_\Gve \|_{L^\infty(\p\GO)} \bigg).
\eeq
It then follows from \eqnref{invA1} and \eqnref{invAy} that
\beq
\psi_\Gve (x)= \frac{4}{E} \left[ \frac{\p \Bu_{0}} {\p \nu}(0)
\sqrt{1-x^2} + \frac{\Gve}{4} \frac{\p^2  \Bu_{0}}  {\p t \p
\nu}(0) x \sqrt{1-x^2} + \Ecal_1(x) \right].
\eeq
Thus we have from \eqnref{psiGve} that
\beq\label{GvfGve}
\Gvf_\Gve (x)= \frac{2}{E} \left[ \frac{\p \Bu_{0}} {\p \nu}(0) \sqrt{\Gve^2 - 4
x^2} + \frac{1}{2} \frac{\p^2  \Bu_{0}}  {\p t \p \nu}(0) x
\sqrt{\Gve^2 - 4 x^2} + \Ecal(x) \right], \quad (x,0) \in
\Gg_\Gve,
\eeq
where $\Ecal(x)= \frac{\Gve}{2} \Ecal_1(\frac{2}{\Gve} x)$ satisfies
\beq\label{Eest}
\| \Ecal \|_{L^\infty(\Gg_\Gve)} \le C \Gve \bigg(\Gve^2 + \| \Bw_\Gve
\|_{L^\infty(\p\GO)} \bigg).
\eeq

Substituting \eqnref{GvfGve} into \eqnref{Neumann_repr} we obtain
\begin{align}
\Bw_\Gve(\Bx) &= \frac{2}{E} \int_{\Gg_\Gve} \frac{\p} {\p \nu_\By} \BN(\Bx, (y,0)) \sqrt{\Gve^2 - 4 y^2} \,dy \frac{\p \Bu_{0}} {\p \nu}(0) \nonumber \\
& \quad + \frac{1}{E} \int_{\Gg_\Gve} \frac{\p} {\p \nu_\By} \BN(\Bx, (y,0)) y \sqrt{\Gve^2 - 4 y^2} \,dy \frac{\p^2  \Bu_{0}}  {\p t \p \nu}(0) \nonumber \\
& \quad + \frac{2}{E} \int_{\Gg_\Gve} \frac{\p} {\p \nu_\By} \BN(\Bx, (y,0)) \Ecal(y) \,dy := I+ II+III , \quad \Bx \in \p \GO . \label{BuGveBx}
\end{align}
Since
$$
\frac{\p}{\p \nu_\By} \BN(\Bx, (y,0))= \frac{\p}{\p \nu_\By} \BN(\Bx, 0) + \frac{\p^2}{\p t_\By \p \nu_\By} \BN(\Bx, 0) y + O(y^2),
$$
we have
\begin{align*}
\int_{\Gg_\Gve} \frac{\p} {\p \nu_\By} \BN(\Bx, (y,0)) \sqrt{\Gve^2 - 4 y^2} \,dy
&= \frac{\p}{\p \nu_\By} \BN(\Bx, 0) \int_{\Gg_\Gve} \sqrt{\Gve^2 - 4 y^2} \,dy \\
& \quad + \frac{\p^2}{\p t_\By \p \nu_\By} \BN(\Bx, 0) \int_{\Gg_\Gve} y \sqrt{\Gve^2 - 4 y^2} \,dy + O(\Gve^4) \\
&= \frac{\pi\Gve^2}{2} \frac{\p}{\p \nu_\By} \BN(\Bx, 0)  + O(\Gve^4),
\end{align*}
and hence \beq I= \frac{\pi\Gve^2}{E} \frac{\p}{\p \nu_\By}
\BN(\Bx, 0) \frac{\p \Bu_{0}} {\p \nu}(0) + O(\Gve^4). \eeq Here
and throughout this paper, $O(\Gve^4)$ is in the sense of the
uniform norm on $\p\GO$. Similarly, one can show that \beq II=
O(\Gve^4). \eeq So we obtain that  \beq\label{BwGve} \Bw_\Gve(\Bx)
= \frac{\pi\Gve^2}{E} \frac{\p}{\p \nu_\By} \BN(\Bx, 0) \frac{\p
\Bu_{0}} {\p \nu}(0) + O(\Gve^4) + III. \eeq In particular, we
have
$$
\| \Bw_\Gve \|_{L^\infty(\p\GO)} \le C (\Gve^2+ |III|).
$$
But, because of \eqnref{Eest}, we arrive at
$$
|III| \le C \Gve \| \Ecal \|_{L^\infty(\Gg_\Gve)} \le C \Gve^2 (\Gve^2 + \| \Bw_\Gve \|_{L^\infty(\p\GO)}),
$$
and hence
$$
\| \Bw_\Gve \|_{L^\infty(\p\GO)} \le C \Gve^2 (1 + \| \Bw_\Gve \|_{L^\infty(\p\GO)}).
$$
So, if $\Gve$ is small enough, then
\beq
\| \Bw_\Gve \|_{L^\infty(\p\GO)} \le C \Gve^2.
\eeq
It then follows from \eqnref{BwGve} that
\beq\label{BwGve2}
\Bw_\Gve(\Bx) = \frac{\pi\Gve^2}{E} \frac{\p}{\p \nu_\By} \BN(\Bx, 0) \frac{\p \Bu_{0}} {\p \nu}(0) + O(\Gve^4) .
\eeq

We obtain the following theorem.
\begin{theorem} \label{thm-neuman} 
Suppose that $\gamma_\Gve$ is a linear crack of size $\Gve$ and $\Bz$ is the center of $\Gg_\Gve$. Then the solution to \eqnref{prob1} has the following asymptotic expansion:
\begin{equation}\label{exp-u}
  (\Bu_{\Gve}- \Bu_0)(\Bx)
  = \frac{\pi\Gve^2}{E} \frac{\p \BN} {\p \nu_\By} (\Bx,\By)\Big|_{\By=\Bz}
\frac{\p \Bu_0}{\p \nu}(\Bz)
  + O(\Gve^{4})
     \end{equation}
uniformly on $\Bx\in \p \Omega$. Here $E$ is the Young's modulus.
\end{theorem}

It is worth emphasizing that in \eqnref{exp-u} the error is $O(\Gve^4)$ and
the $\Gve^3$-term vanishes. One can see from the derivation of \eqnref{exp-u} that
the $\Gve^3$-term vanishes because $\Gg_\Gve$ is a line segment. If it is a curve,
then we expect that the $\Gve^3$-term does not vanish. We also emphasize that \eqnref{exp-u} is a point-wise asymptotic formula, and it can be used to design algorithms to reconstruct cracks from boundary measurements. We can also integrate this formula against the traction $\Bg$ to obtain the asymptotic formula for the perturbation of the elastic energy as we do in the next section.

Similarly, if we consider the Dirichlet problem
 \begin{equation}\label{prob21}
\begin{cases}
\nabla \cdot \Gs(\Bu_{\Gve}) = 0 \quad  & \mbox{in}\; \Omega \setminus \bar{\gamma}_{\Gve}, \\
\Bu_{\Gve} = \Bf   & \mbox{on} \; \p \Omega,\\
\Gs(\Bu_{\Gve})\,\Be^\perp = 0   & \mbox{on} \;
\gamma_{\Gve},
\end{cases}
\end{equation}
and denote the Green function of Lam\'{e} system in $\Omega$ by
$\BG$, then we get the following asymptotic expansion of its
solution $\Bu_{\Gve}$.

\begin{theorem} \label{thm-dirich} 
Suppose that $\gamma_\Gve$ is a linear crack of size $\Gve$, located at $\Bz$. Then the solution to \eqnref{prob21} has the following asymptotic expansion:
\begin{equation}\label{exp-Du}
  \pd{}{\nu}(\Bu_{\Gve}- \Bu_0)(\Bx)
  =  \frac{\pi\Gve^2}{E}
  \frac{\p^2  \BG}{\p \nu_\Bx \p \nu_\By}(\Bx,\By)\Big|_{\By=\Bz}
  \frac{\p \Bu_0}{\p \nu}(\Bz)
   + O(\Gve^{4})
     \end{equation}
uniformly on $\Bx\in \p \Omega$.
\end{theorem}

\section{Topological derivative of the potential energy} \label{sect5}

The elastic potential energy functional of the cracked body is given by \eqnref{energy}, while without the crack the energy functional is given by
\beq
J[\Bu_{0}]
= -\frac{1}{2}\int_{\Omega} \sigma(\Bu_{0}) :  \nabla^s \Bu_{0}.
\eeq
By the divergence theorem we have
\beq
J_{\Gve}[\Bu_{\Gve}] - J[\Bu_{0}] = - \frac{1}{2} \int_{\p\GO}  (\Bu_{\Gve}-\Bu_0) \cdot \Bg \, d\Gs.
\eeq
Thus we obtain from \eqref{exp-u}
 \begin{align*}
 J_{\Gve}[\Bu_{\Gve}] - J[\Bu_{0}] = - \frac{\pi\Gve^2}{2E}
    \frac{\p \Bu_0}{\p \nu}(\Bz)\frac{\p }{\p \nu_\By}
\int_{\p \Omega} \BN (\Bx,\By)
    \Bg(\Bx) \,d\sigma(\Bx)\Big|_{\By=\Bz}  + O(\Gve^{4}).
  \end{align*}
Since
$$
\Bu_0(\By)= \int_{\p \Omega} \BN (\Bx,\By)
    \Bg(\Bx) \,d\sigma(\Bx), \quad \By \in \GO,
$$
we have
\beq\label{energyform}
 J_{\Gve}[\Bu_{\Gve}] - J[\Bu_{0}] = - \frac{\pi\Gve^2}{2E}
\Big|\frac{\p \Bu_0}{\p \nu}(\Bz)\Big|^2 + O(\Gve^{4}).
\eeq

We may write \eqnref{energyform} in terms of the stress intensity factors. The (normalized) stress intensity factors $K_I$ and $K_{II}$ are defined by
\begin{equation}
K_{I}(\Bu_0,\Be):= \sigma(\Bu_0)\Be^{\perp}\cdot \Be^{\perp}\quad \mbox{and}\quad
K_{II}(\Bu_0,\Be):= \sigma(\Bu_0)\Be^\perp \cdot \Be.
\end{equation}
So, we have
\begin{equation}
\sigma (\Bu_0)\Be^{\perp}= K_{I}\Be^{\perp}+K_{II}\Be,
\end{equation}
and hence
\begin{equation}
\Big|\frac{\p \Bu_0}{\p \nu}(\Bz)\Big|^2 = |\sigma (\Bu_0)\Be^{\perp}|^2 = K_{I}^2 + K_{II}^2.
\end{equation}
We obtain the following result.
\begin{theorem}
We have
\beq\label{energyform2}
 J_{\Gve}[\Bu_{\Gve}] - J[\Bu_{0}] = - \frac{\pi\Gve^2}{2E}
(K_{I}^2 + K_{II}^2) + O(\Gve^{4})
\eeq
as $\Gve \to 0$.
\end{theorem}

The topological derivative $D_T J_{\Gve}(\Bz)$ of the potential
energy is defined by \cite{SZ99, SZ01} 
\beq\label{topd} 
D_T J_{\Gve}(\Bz) :=\lim_{\Gve \rightarrow
0}\left(\frac{1}{\rho'(\Gve)} \frac{d}{d\Gve}J_{\Gve}\right), 
\eeq
where $\rho(\Gve)= \pi\Gve^2$. So, one can immediately see from
\eqnref{energyform2} that 
\beq\label{plainstress}
D_T J_{\Gve}(\Bz) = - \frac{1}{2E} (K_{I}^2 + K_{II}^2). 
\eeq 
This formula is in accordance with the one obtained by Novotny {\it et al} in \cite{NDTP} (see also \cite{GN}).
In fact, in those papers the plane strain and the plain stress problems are considered, and  \eqnref{plainstress} is the formula for the latter problem.

\section*{Acknowledgement} Authors would like to thank Andr\'e Novotny for helpful comments on this paper.

\appendix 
\section{Derivation of \eqnref{repDij}} \label{appendix}

The Kelvin matrix \eqref{lamefund} can be rewritten as
\begin{equation}\label{lamefund2}
   \Phi_{ij}(\Bx-\By)=\lambda' \delta_{ij}\log|\Bx-\By|
   + \mu'(x_i-y_i)\frac{\p \log|\Bx-\By| }{\p y_j},\quad
   i,j=1,2,
\end{equation}
where
$$
\lambda' = \frac{\lambda+3\mu}{4\pi \mu(\lambda+2\mu)}, \quad
\mu'= \frac{\lambda+\mu}{4\pi\mu(\lambda+2\mu)}.
$$
Using the operator $T(\p)$ defined by \eqnref{op-T} one can see that
\begin{equation}\label{derfund-rule}
 \frac{\p \Phi}{\p \nu_\By}(\Bx-\By)
= (T (\p_\By)\Phi(\Bx-\By))^T,
\end{equation}
or
\begin{equation}\label{defopr}
\displaystyle\Big(\frac{\p \Phi}{\p \nu_\By}(\Bx-\By)\Big)_{kj}:=
\sum_{l=1}^2 T_{jl}(\p_\By)\Phi_{lk}(\Bx-\By), \quad k,j=1,2.
\end{equation}
We use the formulas
\begin{eqnarray*}
 \frac{\p}{\p y_i}\log|\Bx-\By| &=&-\frac{x_i-y_i}{|\Bx-\By|^2},\\
\frac{\p^2}{\p y_i^2}\log|\Bx-\By| &=&-2\frac{(x_i-y_i)^2}{|\Bx-\By|^4}+ \frac{1}{|\Bx-\By|^2},\\
\frac{\p^2}{\p y_i\p y_j}\log|\Bx-\By| &=&-2\frac{(x_i-y_i)(x_j-y_j)}{|\Bx-\By|^4} \quad\mbox{if } i \neq j.
  \end{eqnarray*}
By \eqref{defopr}, we have
\begin{eqnarray*}
\Big(\frac{\p \Phi}{\p \nu_\By}\Big)_{11}
&=& \left((\lambda +2\mu) n_1 \frac{\p }{\p {y_1}} + \mu n_2\frac{\p }{\p {y_2}}\right)\left(\lambda' \log|\Bx-\By| + \mu'(x_1-y_1)\frac{\p \log|\Bx-\By|}{\p {y_1}}\right)\\
&& + \left( \mu n_2 \frac{\p }{\p {y_1}} + \lambda n_1 \frac{\p }{\p {y_2}}\right)\mu'(x_1-y_1)\frac{\p \log|\Bx-\By|}{\p {y_2}}\\
&=& \lambda'(\lambda +2\mu) n_1 \frac{\p \log|\Bx-\By|}{\p {y_1}} + \lambda'\mu n_2\frac{\p \log|\Bx-\By|}{\p {y_2}}\\
&& +\,  \mu'(\lambda +2\mu) n_1 \left(-\frac{\p \log|\Bx-\By|}{\p {y_1}}+(x_1-y_1)\frac{\p^2 \log|\Bx-\By|}{\p {y_1}^2}\right)\\
&& - \mu\mu'n_2\frac{\p \log|\Bx-\By|}{\p {y_2}} + 2\mu\mu'n_2(x_1-y_1)\frac{\p^2 \log|\Bx-\By|}{\p {y_1}\p{y_2}} \\ &&+\lambda\mu'(x_1-y_1)\frac{\p^2 \log|\Bx-\By|}{\p {y_2}^2}.
  \end{eqnarray*}
Since $\Delta \log|\Bx-\By|=0$ for $\Bx \neq \By$, we have
\begin{eqnarray*}
\Big(\frac{\p \Phi}{\p \nu_\By}\Big)_{11}
 &=& (\lambda +2\mu)(\lambda'-\mu') n_1 \frac{\p \log|\Bx-\By| }{\p y_1}
+\mu (\lambda'-\mu') n_2\frac{\p\log|\Bx-\By| }{\p y_2} \\
&& +2\mu\mu' \left( n_1(x_1-y_1)\frac{\p^2 \log|\Bx-\By|}{\p {y_1}^2}+n_2(x_1-y_1)\frac{\p^2 \log|\Bx-\By|}{\p {y_1}\p{y_2}}\right).
 \end{eqnarray*}
Since
$$
(\lambda +2\mu)(\mu'-\lambda')+2\mu\mu'=-\frac{\mu}{2\pi (\lambda +2\mu) }=\mu(\mu'-\lambda'),
$$
we obtain
\begin{eqnarray*}
\Big(\frac{\p \Phi}{\p \nu_\By}\Big)_{11}
=\left[ \mu(\mu'-\lambda') - 4\mu\mu'\frac{(x_1-y_1)^2}{|\Bx-\By|^2}\right]
 \sum_{l=1}^2 n_l\frac{x_l-y_l}{|\Bx-\By|^2}.
 \end{eqnarray*}

Similarly, we can compute
\begin{eqnarray*}
\Big(\frac{\p \Phi}{\p \nu_\By}\Big)_{12}
&=& \left(\lambda n_2 \frac{\p }{\p {y_1}} + \mu n_1\frac{\p }{\p {y_2}}\right)\left(\lambda' \log|\Bx-\By|  + \mu'(x_1-y_1)\frac{\p \log|\Bx-\By|}{\p {y_1}}\right)\\
&& + \left(\mu n_1 \frac{\p }{\p {y_1}}
 + (\lambda +2\mu) n_2 \frac{\p }{\p {y_2}}\right)\mu'(x_1-y_1)\frac{\p \log|\Bx-\By|}{\p {y_2}}\\
&=&\lambda(\lambda' - \mu')n_2 \frac{\p \log|\Bx-\By|}{\p {y_1}}
+\mu( \lambda' -\mu')n_1\frac{\p \log|\Bx-\By|}{\p {y_2}}\\
&&+2\mu\mu'n_1(x_1-y_1) \frac{\p ^2\log|\Bx-\By|}{\p {y_1} \p{y_2}}
+ 2\mu\mu' n_2 (x_1-y_1)\frac{\p ^2\log|\Bx-\By|}{\p {y_2}^2}.
  \end{eqnarray*}
Since $\lambda(\mu'-\lambda')+2\mu\mu'=\mu(\lambda'-\mu')$, we have
\begin{eqnarray*}
\Big(\frac{\p \Phi}{\p \nu_\By}\Big)_{12}
&=&[\lambda(\mu'-\lambda')+2\mu\mu']n_2 \frac{x_1-y_1}{|\Bx-\By|^2}
+\mu(\mu'- \lambda')n_1\frac{x_2-y_2}{|\Bx-\By|^2} \\
&&-4\mu\mu'n_1 \frac{(x_1-y_1)^2(x_2-y_2)}{|\Bx-\By|^4}
-4\mu\mu' n_2 \frac{(x_1-y_1)(x_2-y_2)^2}{|\Bx-\By|^4}\\
 &=& -4\mu\mu'\frac{(x_1-y_1)(x_2-y_2)}{|\Bx-\By|^2}
 \sum_{l=1}^2 n_l\frac{x_l-y_l}{|\Bx-\By|^2} \\
&&  \quad -\,\mu(\mu'-\lambda') \left[\frac{n_2(x_1-y_1) - n_1(x_2-y_2)}{|\Bx-\By|^2}\right].
  \end{eqnarray*}

We also have
\begin{align*}
\Big(\frac{\p \Phi}{\p \nu_\By}\Big)_{22}
& = \left(\lambda n_2 \frac{\p }{\p {y_1}} + \mu n_1\frac{\p }{\p {y_2}}\right)\mu'(x_2-y_2)\frac{\p \log|\Bx-\By|}{\p {y_1}}\\
&+\left(\mu n_1\frac{\p }{\p {y_1}}+(\lambda +2\mu)n_2\frac{\p }{\p {y_2}}\right) \left(\lambda' \log|\Bx-\By|  + \mu'(x_2-y_2)\frac{\p \log|\Bx-\By|}{\p {y_2}}\right)\\
& = \left[ \mu(\mu'-\lambda') - 4\mu\mu'\frac{(x_2-y_2)^2}{|\Bx-\By|^2}\right]
 \sum_{l=1}^2 n_l\frac{x_l-y_l}{|\Bx-\By|^2}.
 \end{align*}
 This proves \eqref{repDij}.


\end{document}